\documentclass[11pt]{article}
\usepackage{amsmath, amssymb, amsfonts}
\usepackage{epsfig}

\def\be{\begin{equation}}
\def\ee{\end{equation}}
\def\bea{\begin{eqnarray}}
\def\eea{\end{eqnarray}}
\def\bes{\begin{eqnarray*}}
\def\ees{\end{eqnarray*}}

\def\nn{\nonumber}
\def\<{\langle}
\def\>{\rangle}
\def\lb{\label}
\def\bs{\setminus}

\def\R{{\mathbb{R}}}
\def\C{{\bf C}}
\def\Z{{\mathbb{Z}}}
\def\N{{\mathbb{N}}}
\def\U{{\bf U}}
\def\Q{{\mathbb{Q}}}

\def\RP{{\mathbb{R}P^{n}}}

\def\aa{{\alpha}}
\def\bb{{\beta}}
\def\ga{{\gamma}}

\def\th{{\theta}}
\def\om{{\omega}}
\def\Om{{\Omega}}
\def\ep{{\epsilon}}
\def\lm{{\lambda}}
\def\Lm{{\Lambda}}

\def\sg{{\sigma}}

\def\vf{{\varphi}}

\def\K{{\cal K}}
\def\P{{\cal P}}

\def\ker{{\rm ker}}

\def\hv{{\rm hv}}
\def\diag{{\rm diag}}

\def\rank{{\rm rank}}
\def\Sp{{\rm Sp}}
\def\mod{{\rm mod}}
\def\CG{{\rm CG}}

\def\dm{{\diamond}}
\def\ol#1{\overline{#1}}  

\title{The optimal lower bound estimation of the number of closed geodesics on Finsler compact space form $S^{2n+1}/ \Gamma$}
\author{Hui Liu\thanks{Partially supported by NSFC (Nos. 11401555, 11771341),
Anhui Provincial Natural Science Foundation (No. 1608085QA01).
E-mail: huiliu00031514@whu.edu.cn.} \\\\
School of Mathematics and Statistics, Wuhan University,
\\Wuhan 430072, Hubei, China\\}
\date{}
\begin{document}
\newtheorem{definition}{Definition}[section]
\newtheorem{theorem}{Theorem}[section]
\newtheorem{lemma}{Lemma}[section]
\newtheorem{corollary}{Corollary}[section]
\newtheorem{example}{Example}[section]
\newtheorem{property}{Property}[section]
\newtheorem{proposition}{Proposition}[section]
\newtheorem{remark}{Remark}[section]

\newcommand{\qed}{\nolinebreak\hfill\rule{2mm}{2mm}
\par\medbreak}
\newcommand{\Proof}{\par\medbreak\it Proof: \rm}
\newcommand{\rem}{\par\medbreak\it Remark: \rm}
\newcommand{\defi}{\par\medbreak\it Definition : \rm}
\renewcommand{\thefootnote}{\arabic{footnote}}

\maketitle

\begin{abstract}
{\it Let $M=S^{2n+1}/ \Gamma$, $\Gamma$ is a finite group which acts
freely and isometrically on the $(2n+1)$-sphere and therefore $M$ is diffeomorphic to a compact space form.
In this paper, we first investigate Katok's famous example about irreversible Finsler metrics on the spheres
to study the topological structure of the contractible component of the free loop space on
the compact space form $M$, then we apply the result to establish the resonance identity for
homologically visible contractible minimal closed geodesics on every Finsler compact space form $(M, F)$
when there exist only finitely many distinct contractible minimal closed geodesics on $(M, F)$.
As its applications, using this identity and the enhanced common index jump
theorem for symplectic paths proved by Duan, Long and Wang in \cite{DLW2}, we show that there exist at least
$2n+2$ distinct closed geodesics on every compact space form $S^{2n+1}/ \Gamma$ with
a bumpy irreversible Finsler metric $F$ under some natural curvature condition, which is the optimal lower
bound due to Katok's example.  }
\end{abstract}



\renewcommand{\theequation}{\thesection.\arabic{equation}}
\makeatletter\@addtoreset{equation}{section}\makeatother

\setcounter{equation}{0}
\section{Introduction}
Let $M=S^{2n+1}/ \Gamma$, $\Gamma$ is a finite group which acts freely and
isometrically on the $(2n+1)$-sphere and therefore $M$ is diffeomorphic to a compact space form
which is typically a non-simply connected manifold. In papers \cite{Tai2016} and \cite{LLX},
the existence of at least two distinct non-contractible closed geodesics on every bumpy $S^n/ \Gamma$ with $n\geq2$
was proved, however, this paper is concerned with the total number of closed geodesics on Finsler $S^{2n+1}/ \Gamma$,
the main ingredients are the investigations of the topological
structure of the contractible component of the free loop space
on $S^{2n+1}/ \Gamma$ which yields a new resonance identity for
homologically visible contractible minimal closed geodesics on Finsler $S^{2n+1}/ \Gamma$
, and the enhanced common index jump theorem for symplectic paths discovered by Duan, Long and Wang in \cite{DLW2}.

A closed curve on a Finsler manifold is a closed geodesic if it is
locally the shortest path connecting any two nearby points on this
curve. As usual, on any Finsler manifold
$(M, F)$, a closed geodesic $c:S^1=\R/\Z\to M$ is {\it prime}
if it is not a multiple covering (i.e., iteration) of any other closed geodesics.
 Here the $m$-th iteration $c^m$ of $c$ is defined
by $c^m(t)=c(mt)$.  The inverse curve $c^{-1}$ of $c$ is defined by
$c^{-1}(t)=c(1-t)$ for $t\in \R$.  Note that unlike Riemannian manifold,
the inverse curve $c^{-1}$ of a closed geodesic $c$
on a irreversible Finsler manifold need not be a geodesic.
We call two prime closed geodesics
$c$ and $d$ {\it distinct} if there is no $\th\in (0,1)$ such that
$c(t)=d(t+\th)$ for all $t\in\R$.
For a closed geodesic $c$ on $(M,\,F)$, denote by $P_c$
the linearized Poincar\'{e} map of $c$. Recall that a Finsler metric $F$ is {\it bumpy} if all the closed geodesics
on $(M, \,F)$ are non-degenerate, i.e., $1\notin \sigma(P_c)$ for any closed
geodesic $c$.

Let $\Lm M$ be the free loop space on $M$ defined by
\begin{equation*}
\label{LambdaM}
  \Lambda M=\left\{\gamma: S^{1}\to M\mid \gamma\ {\rm is\ absolutely\ continuous\ and}\
                        \int_{0}^{1}F(\gamma,\dot{\gamma})^{2}dt<+\infty\right\},
\end{equation*}
endowed with a natural structure of Riemannian Hilbert manifold on which the group $S^1=\R/\Z$ acts continuously by
isometries (cf. Shen \cite{Shen2001}).

It is well known (cf. Chapter 1 of Klingenberg \cite{Kli1978}) that $c$ is a closed geodesic or a constant curve
on $(M,F)$ if and only if $c$ is a critical point of the energy functional
\begin{equation*}
\label{energy}
E(\gamma)=\frac{1}{2}\int_{0}^{1}F(\gamma,\dot{\gamma})^{2}dt.
\end{equation*}
Based on it, many important results on this subject have been obtained
(cf. \cite{Ano},  \cite{Hin1984}-\cite{Hingston1993},
\cite{Rad1989}-\cite{Rad1992}). In particular, in 1969 Gromoll and Meyer \cite{GM1969JDG}
used Morse theory and Bott's index iteration formulae \cite{Bott1956}
to establish the existence of infinitely many distinct closed geodesics on $M$, when the Betti number sequence
$\{{\beta}_k(\Lm M;\mathbb{Q})\}_{k\in\Z}$ is unbounded. Then Vigu$\acute{e}$-Poirrier and Sullivan \cite{VS1976}
further proved in 1976 that for a compact simply connected manifold $M$, the Gromoll-Meyer condition holds if and
only if $H^{*}(M;\mathbb{Q})$ is generated by more than one element.

However, when $\{{\beta}_k(\Lm M;\mathbb{Q})\}_{k\in\Z}$ is bounded, the problem is quite complicated.
In 1973, Katok \cite{Katok1973} endowed some irreversible Finsler metrics to the compact rank one symmetric spaces
$$  S^{n},\ \mathbb{R}P^{n},\ \mathbb{C}P^{n},\ \mathbb{H}P^{n}\ \text{and}\ {\rm CaP}^{2},  \label{mflds}$$
each of which possesses only finitely many distinct prime closed geodesics
(cf. also Ziller \cite{Ziller1977},\cite{Ziller1982}).
On the other hand, Franks \cite{Franks1992} and Bangert \cite{Bangert1993} together proved that
there are always infinitely many distinct closed geodesics on every Riemannian sphere $S^2$
(cf. also  Hingston \cite{Hingston1993}, Klingenberg \cite{Kli1995}).
These results imply that the metrics play an important role on the multiplicity of closed geodesics on those manifolds.

In 2004,  Bangert and Long \cite{BL2010} (published in 2010) proved the existence of at least two distinct
closed geodesics on every Finsler $S^2$. Subsequently, such a multiplicity result for
$S^{n}$ with a bumpy Finsler metric  was proved by Duan and Long
\cite{DL2007} and Rademacher \cite{Rad2010} independently. Furthermore in a recent paper \cite{DLW1},
Duan, Long and Wang  proved the same conclusion for any compact simply-connected bumpy Finsler manifold,
and in \cite{DLW2} they obtained the optimal lower bound estimation of the number of closed geodesics on
any compact simply-connected bumpy Finsler manifold under some curvature conditions.
We refer the readers to \cite{DL2010}, \cite{HiR}, \cite{LD2009}, \cite{Rad2007} \cite{Wang2008}-\cite{Wang2012}
and the references therein for more interesting results and the survey papers of Long \cite{lo2006},
Taimanov \cite{Tai2010}, Burns and Matveev \cite{BM2013} and Oancea \cite{Oancea2014}
for more recent progresses on this subject.

As for the multiplicity of closed geodesics on non-simply connected manifolds whose
free loop space possesses bounded Betti number sequence, Ballman et al. \cite{BTZ1} proved in 1981
that every Riemannian manifold with the fundamental group being a nontrivial finitely cyclic
group and possessing a generic metric has infinitely many distinct closed geodesics. In 1984,
Bangert and Hingston \cite{BaH} proved that any Riemannian manifold with the fundamental group being an
infinite cyclic group has infinitely many distinct closed geodesics. Since then, there seem to be
very few works on the multiplicity of closed geodesics on non-simply connected manifolds.
The main reason is that the topological structures of the free loop spaces on these manifolds
are not well known, so that the classical Morse theory is difficult to be applicable.

Motivated by the studies on simply connected manifolds, in particular,
the resonance identity proved by Rademacher \cite{Rad1989},
and based on Westerland's works \cite{West2005}, \cite{West2007}
on loop homology of $\RP$, Xiao and Long \cite{XL2015}
in 2015 investigated the topological structure of the non-contractible
loop space and established the resonance identity for the non-contractible
closed geodesics on $\R P^{2n+1}$ by use of
$\Z_2$ coefficient homology. As an application, Duan, Long and Xiao \cite{DLX2015}
proved the existence of at least two
distinct non-contractible closed geodesics on $\R P^{3}$ endowed with a bumpy
and irreversible Finsler metric. Subsequently in \cite{Tai2016}, Taimanov
used a quite different method from \cite{XL2015} to compute the rational equivariant cohomology of
 the non-contractible loop spaces in compact space forms $S^n/ \Gamma$ and
proved the existence of at least two distinct non-contractible closed geodesics on $\mathbb{R}P^2$
endowed with a bumpy irreversible Finsler metric.
Then in \cite{Liu}, Liu combined Fadell-Rabinowitz index theory with Taimanov's topological results
to get many multiplicity
results of non-contractible closed geodesics on positively curved Finsler $\RP$.
In \cite{LX},  Liu and Xiao established
the resonance identity for the non-contractible closed geodesics on $\mathbb{R}P^n$, and
together with \cite{DLX2015} and \cite{Tai2016} proved the existence of at least two distinct
non-contractible closed geodesics on every bumpy $\mathbb{R}P^n$ with $n\geq2$.
Furthermore, Liu, Long and Xiao \cite{LLX}
proved that every bumpy Finsler compact space form $S^n/ \Gamma$ possesses two distinct
closed geodesics in each of its nontrivial classes.

Based on the works of \cite{DLW2} and \cite{LLX}, it is natural to ask whether we can improve the lower bound
of the total number of closed geodesics on bumpy Finsler compact space forms under some natural conditions.
Note that the only non-trivial group which acts freely on $S^{2n}$ is $\Z_2$
and $S^{2n}/ \Z_2=\mathbb{R}P^{2n}$(cf. P.5 of \cite{Tai2016}).
In \cite{GGM}, Ginzburg, Gurel and Macarini have obtained the optimal lower bound
of the total number of closed geodesics on bumpy Finsler  $\mathbb{R}P^{2n}$,
so the left interesting part about this problem is to estimate the lower bound
of the total number of closed geodesics on bumpy Finsler compact space form $S^{2n+1}/ \Gamma$.
This paper is devoted to answer this question. To this end, we first investigate Katok's example
about irreversible Finsler metrics on the spheres
and  compute the $S^1$-equivariant Betti number sequence of the contractible component of the free loop space on
the compact space form $S^{2n+1}/ \Gamma$ in Section 4, and then in Section 5 we use this topological result to
establish the following resonance identity. Note that a contractible closed geodesic $c$ is called {\it minimal} if it is not an iteration of any other contractible closed geodesics.
\begin{theorem}\label{Thm1.1} Let $M=S^{2n+1}/ \Gamma$ and $e$ be the identity in $\pi_1(M)$.
Suppose the Finsler manifold $(M,F)$ possesses only finitely many distinct contractible minimal closed geodesics,
among which we denote the distinct homologically visible contractible minimal
closed geodesics by $c_1, \ldots, c_r$  for some integer $r>0$. Then we have
\bea  \sum_{j=1}^{r}\frac{\hat{\chi}(c_j)}{\hat{i}(c_j)} = \bar{B}(\Lambda_e M, \Lambda^0 M;\Q) =
    \frac{n+1}{2n},     \lb{reident1}\eea
where $\Lm_{e} M$ is the contractible loop space of $M$,
$\Lambda^0 M  =\{{\rm constant\;point\;curves\;in\;}M\}\cong M$
and the mean Euler number $\hat{\chi}(c_j)$ of $c_j$ is defined by
$$  \hat{\chi}(c_j) = \frac{1}{n_j}\sum_{m=1}^{n_j}\sum_{l=0}^{4n}(-1)^{l+i(c_{j}^{m})}
k_{l}^{\epsilon(c_j^{m})}(c_{j}^{m})\in\Q, $$
and $n_j=n_{c_j}$ is the analytical period of $c_j$, $k_{l}^{\epsilon(c_j^{m})}(c_{j}^m)$
is the local homological type number of $c_{j}^{m}$,
$i(c_{j})$ and $\hat{i}(c_j)$ are the Morse index and mean index of $c_j$ respectively.

In particular, if the Finsler metric $F$ on $M=S^{2n+1}/ \Gamma$ is bumpy,
then (\ref{reident1}) has the following simple form
\bea  \sum_{j=1}^{r}\left(\frac{(-1)^{i(c_{j})}k_0^{\epsilon(c_j)}(c_{j})+(-1)^{i(c_{j}^2)}
k_0^{\epsilon(c_j^2)}(c_{j}^2)}{2}\right)\frac{1}{\hat{i}(c_{j})}
=\frac{n+1}{2n}.  \lb{breident1}\eea
\end{theorem}

Based on Theorem \ref{Thm1.1},  we use Morse theory and the enhanced common index jump
theorem for symplectic paths proved by Duan, Long and Wang in \cite{DLW2} to
prove our main multiplicity results of contractible closed geodesics on $(S^{2n+1}/ \Gamma, F)$.
\begin{theorem}
\label{mainresult}
Let $M=(S^{2n+1}/ \Gamma, F)$ be a bumpy Finsler compact space form. If the number
of contractible closed geodesics is finite and the flag curvature $K$ satisfies
$\left(\frac{\lambda}{1+\lambda}\right)^2<K\le 1$,
then there exist at least $(2n+2)$ distinct contractible closed geodesics with even Morse indices, and $2n$ ones
of them are non-hyperbolic.
\end{theorem}
\begin{remark}
Due to Katok's famous example, the lower bound of the number of closed geodesics obtained
in the above theorem is optimal, cf. \cite{Ziller1982} and Section 4 below.
\end{remark}

This paper is organized as follows. In Section 2, we review the critical point theory of closed  geodesics
and apply the splitting theorem on critical modules to the contractible component of the free loop space of
$(S^{2n+1}/ \Gamma, F)$. Then in Section 3, we recall the precise iteration formulae
and the enhanced common index jump theorem of symplectic paths,
which also work for Morse indices of orientable closed geodesics.
In Section 4, we investigate Katok's famous example about irreversible Finsler metrics on the spheres
to study the topological structure of the contractible component of the free loop space on
the compact space form $S^{2n+1}/ \Gamma$, then in Section 5 we apply the result to establish
the resonance identity in Theorem \ref{Thm1.1}. Finally in Section 6,
we give the proof of our main Theorem \ref{mainresult}.

In this paper, let $\N$, $\N_0$,  $\Z$, $\Q$ and $\Q^{c}$ denote the sets of natural integers, non-negative
integers, integers, rational numbers and  irrational numbers respectively.
We also use  notations $E(a)=\min\{k\in\Z\,|\,k\ge a\}$,
$[a]=\max\{k\in\Z\,|\,k\le a\}$, $\varphi(a)=E(a)-[a]$ and $\{a\}=a-[a]$ for any $a\in\R$.
Throughout this paper, we use $\Q$ coefficients for all homological and cohomological modules.

\setcounter{equation}{0}
\section{Critical point theory for closed geodesics }
Let $M=(M,F)$ be a compact Finsler manifold, the space
$\Lambda=\Lambda M$ of $H^1$-maps $\gamma:S^1\rightarrow M$ has a
natural structure of Riemannian Hilbert manifolds on which the
group $S^1=\R/\Z$ acts continuously by isometries. This action is defined by
$(s\cdot\gamma)(t)=\gamma(t+s)$ for all $\gamma\in\Lm$ and $s,
t\in S^1$. For any $\gamma\in\Lambda$, the energy functional is
defined by
\be E(\gamma)=\frac{1}{2}\int_{S^1}F(\gamma(t),\dot{\gamma}(t))^2dt.
\lb{2.1}\ee
It is $C^{1,1}$ and invariant under the $S^1$-action. The
critical points of $E$ of positive energies are precisely the closed geodesics
$\gamma:S^1\to M$. The index form of the functional $E$ is well
defined along any closed geodesic $c$ on $M$, which we denote by
$E''(c)$. As usual, we denote by $i(c)$ and
$\nu(c)$ the Morse index and nullity of $E$ at $c$.
For a closed geodesic $c$ we set $ \Lm(c)=\{\ga\in\Lm\mid E(\ga)<E(c)\}$.

Recall that respectively the mean index $\hat{i}(c)$ and the $S^1$-critical modules of $c^m$ are defined by
\be \hat{i}(c)=\lim_{m\rightarrow\infty}\frac{i(c^m)}{m}, \quad \overline{C}_*(E,c^m)
   = H_*\left((\Lm(c^m)\cup S^1\cdot c^m)/S^1,\Lm(c^m)/S^1; \Q\right).\lb{2.2}\ee

We call a closed geodesic $c$ satisfying the isolation condition, if
the following holds:

{\bf (Iso)  For all $m\in\N$ the orbit $S^1\cdot c^m$ is an
isolated critical orbit of $E$. }

Note that if the number of prime closed geodesics on a Finsler manifold
is finite, then all the closed geodesics satisfy (Iso).

If $c$ has multiplicity $m$, then the subgroup $\Z_m=\{\frac{n}{m}\mid 0\leq n<m\}$
of $S^1$ acts on $\overline{C}_*(E,c)$. As studied in p.59 of \cite{Rad1992},
for all $m\in\N$, let
$H_{\ast}(X,A)^{\pm\Z_m}
   = \{[\xi]\in H_{\ast}(X,A)\,|\,T_{\ast}[\xi]=\pm [\xi]\}$,
where $T$ is a generator of the $\Z_m$-action.
On $S^1$-critical modules of $c^m$, the following lemma holds:
\begin{lemma}
\label{Rad1992} (cf. Satz 6.11 of \cite{Rad1992} and \cite{BL2010})
Suppose $c$ is
a prime closed geodesic on a Finsler manifold $M$ satisfying (Iso). Then
there exist $U_{c^m}^-$ and $N_{c^m}$, the so-called local negative
disk and the local characteristic manifold at $c^m$ respectively,
such that $\nu(c^m)=\dim N_{c^m}$ and
\bea \overline{C}_q( E,c^m)
&\equiv& H_q\left((\Lm(c^m)\cup S^1\cdot c^m)/S^1, \Lm(c^m)/S^1\right)\nn\\
&=& \left(H_{i(c^m)}(U_{c^m}^-\cup\{c^m\},U_{c^m}^-)
    \otimes H_{q-i(c^m)}(N_{c^m}\cup\{c^m\},N_{c^m})\right)^{+\Z_m}, \nn
\eea

(i) When $\nu(c^m)=0$, there holds
\bea \overline{C}_q( E,c^m) = \left\{\begin{array}{ll}
     \Q, &\quad {\it if}\;\; i(c^m)-i(c)\in 2\Z\;\;{\it and}\;\;
                   q=i(c^m),\;  \cr
     0, &\quad {\it otherwise}, \\ \end{array}\right. \nn \eea

(ii) When $\nu(c^m)>0$, there holds
$$ \overline{C}_q( E,c^m)=H_{q-i(c^m)}(N_{c^m}\cup\{c^m\},N_{c^m})^{\ep(c^m)\Z_m}, $$
where $\ep(c^m)=(-1)^{i(c^m)-i(c)}$.
\end{lemma}

In the following, we let $M=S^{2n+1}/ \Gamma$ and $e$ be the identity in $\pi_1(M)$,
$\Gamma$ acts freely and isometrically on the $(2n+1)$-sphere and therefore
$M$ is diffeomorphic to a compact space form. Then the free loop space $\Lambda M$ possesses
a natural decomposition\bea \Lambda M=\bigsqcup_{g\in \pi_1(M)}\Lambda_g M,\nn\eea
where $\Lambda_g M$ is the connected components of $\Lambda M$ whose elements are homotopic
to $g$. We set $\Lambda_{e}(c) = \{\gamma\in \Lambda_{e}M\mid E(\gamma)<E(c)\}$.
Note that for a contractible minimal closed geodesic $c$, $c^m\in\Lambda_e M$
for every $m\in\N$.

Now we restrict the energy functional $E$ on the contractible component $\Lambda_e M$ and study the Morse theory on $\Lambda_e M$.
Thus we define the $S^1$-critical modules of $c^m$ for $E|_{\Lambda_e M}$ as
\be \overline{C}_*(E,c^m; [e])
   = H_*\left((\Lm_e(c^m)\cup S^1\cdot c^m)/S^1,\Lm_e(c^m)/S^1; \Q\right).\nn\ee
Then by the same proof of Lemma 2.1 we have:
\begin{proposition}
\label{Rad1992'}
Suppose $c$ is a contractible minimal closed geodesic on Finsler
$M=S^{2n+1}/ \Gamma$ satisfying (Iso). Then
there exist $U_{c^m}^-$ and $N_{c^m}$, the so-called local negative
disk and the local characteristic manifold at $c^m$ respectively,
such that $\nu(c^m)=\dim N_{c^m}$ and
\bea \overline{C}_q( E,c^m;[e])
&\equiv& H_q\left((\Lm_e(c^m)\cup S^1\cdot c^m)/S^1, \Lm_e(c^m)/S^1\right)\nn\\
&=& \left(H_{i(c^m)}(U_{c^m}^-\cup\{c^m\},U_{c^m}^-)
    \otimes H_{q-i(c^m)}(N_{c^m}\cup\{c^m\},N_{c^m})\right)^{+\Z_m}, \nn
\eea

(i) When $\nu(c^m)=0$, there holds
\bea \overline{C}_q( E,c^m;[e]) = \left\{\begin{array}{ll}
     \Q, &\quad {\it if}\;\; i(c^m)-i(c)\in 2\Z\;\;{\it and}\;\;
                   q=i(c^m),\;  \cr
     0, &\quad {\it otherwise}, \\ \end{array}\right. \nn \eea

(ii) When $\nu(c^m)>0$, there holds
$$ \overline{C}_q( E,c^m; [e])=H_{q-i(c^m)}(N_{c^m}\cup\{c^m\},N_{c^m})^{\ep(c^m)\Z_m}, $$
where $\ep(c^m)=(-1)^{i(c^m)-i(c)}$.
\end{proposition}

\medskip

As usual, for $m\in\N$ and $l\in\Z$ we define the local homological type numbers of $c^m$ by
\be k_{l}^{\epsilon(c^{m})}(c^{m})
= \dim H_{l}(N_{c^m}\cup\{c^m\},N_{c^m})^{\ep(c^m)\Z_m}.  \lb{CGht1}\ee

Based on works of Rademacher in \cite{Rad1989}, Long and Duan in \cite{LD2009} and \cite{DL2010},
we define the {\it analytical period} $n_c$ of the closed geodesic $c$ by
\be n_c = \min\{j\in \N\,|\,\nu(c^j)=\max_{m\ge 1}\nu(c^m),\;i(c^{m+j})-i(c^m)\in2\Z,\;
                  \forall\,m\in \N\}. \lb{CGap1}\ee
Then by \cite{LD2009} and \cite{DL2010}, we have
\be  k_{l}^{\ep(c^{m+hn_c})}(c^{m+hn_c}) = k_{l}^{\ep(c^m)}(c^m), \qquad \forall\;m,\;h\in \N,\;l\in\Z.  \lb{CGap2}\ee
For more detailed properties of the analytical period $n_c$ of a closed geodesic $c$, we refer readers to
the two Section 3s in \cite{LD2009} and \cite{DL2010}.

\section{The enhanced common index jump theorem of symplectic paths}

In \cite{lo1999} of 1999, Y. Long established the basic normal form decomposition of symplectic matrices.
Based on this result he further established the precise iteration formulae of indices of
symplectic paths in \cite{lo2000} of 2000. Note that this index iteration formulae works for Morse indices
of iterated closed geodesics (cf. \cite{LL2002}, \cite{Liu2005} and Chap. 12 of \cite{lo2002}).
Since every closed geodesic on odd dimensional Finsler manifolds is orientable, then by Theorem 1.1 of \cite{Liu2005}
the initial Morse index of a closed geodesic $c$ on a Finsler compact space form $S^{2n+1}/ \Gamma$ coincides with the index of a
corresponding symplectic path.

As in \cite{lo2000}, denote by
\bea
N_1(\lm, b) &=& \left(
  \begin{array}{ll}\lm & b\cr
                                0 & \lm\\ \end{array}
  \right), \qquad {\rm for\;}\lm=\pm 1, \; b\in\R, \lb{3.1}\\
D(\lm) &=& \left(
  \begin{array}{ll}\lm & 0\cr
                      0 & \lm^{-1}\\ \end{array}
  \right), \qquad {\rm for\;}\lm\in\R\bs\{0, \pm 1\}, \lb{3.2}\\
R(\th) &=& \left(
  \begin{array}{ll}\cos\th & -\sin\th \cr
                           \sin\th & \cos\th\\ \end{array}
  \right), \qquad {\rm for\;}\th\in (0,\pi)\cup (\pi,2\pi), \lb{3.3}\eea
\bea N_2(e^{\th\sqrt{-1}}, B) &=& \left(
  \begin{array}{ll} R(\th) & B \cr
                  0 & R(\th)\\ \end{array}
  \right), \qquad {\rm for\;}\th\in (0,\pi)\cup (\pi,2\pi)\;\; {\rm and}\; \nn\\
        && \qquad B=\left(
  \begin{array}{ll}b_1 & b_2\cr
                                  b_3 & b_4\\ \end{array}
  \right)\; {\rm with}\; b_j\in\R, \;\;
                                         {\rm and}\;\; b_2\not= b_3. \lb{3.4}\eea
Here $N_2(e^{\th\sqrt{-1}}, B)$ is non-trivial if $(b_2-b_3)\sin\theta<0$, and trivial
if $(b_2-b_3)\sin\theta>0$.

As in \cite{lo2000}, the $\diamond$-sum (direct sum) of any two real matrices is defined by
\bea \left(
  \begin{array}{ll}A_1 & B_1\cr C_1 & D_1\\ \end{array}
  \right)_{2i\times 2i}\diamond
      \left(
  \begin{array}{ll}A_2 & B_2\cr C_2 & D_2\\ \end{array}
  \right)_{2j\times 2j}
=\left(
  \begin{array}{llll}A_1 & 0 & B_1 & 0 \cr
                                   0 & A_2 & 0& B_2\cr
                                   C_1 & 0 & D_1 & 0 \cr
                                   0 & C_2 & 0 & D_2\\ \end{array}
  \right).\lb{3.5}\eea

For every $M\in\Sp(2n)$, the homotopy set $\Omega(M)$ of $M$ in $\Sp(2n)$ is defined by
$$ \Om(M)=\{N\in\Sp(2n)\,|\,\sg(N)\cap\U=\sg(M)\cap\U\equiv\Gamma\;\mbox{and}
                    \;\nu_{\om}(N)=\nu_{\om}(M)\, \forall\om\in\Gamma\}, $$
where $\sg(M)$ denotes the spectrum of $M$,
$\nu_{\om}(M)\equiv\dim_{\C}\ker_{\C}(M-\om I)$ for $\om\in\U$.
We denote by $\Om^0(M)$ the path connected component of $\Om(M)$ containing $M$.
\begin{lemma}(cf. \cite{lo2000}, Lemma 9.1.5 and List 9.1.12 of \cite{lo2002})
For $M\in\Sp(2n)$ and $\om\in\U$, the splitting number $S_M^\pm(\om)$
(cf. Definition 9.1.4 of \cite{lo2002}) satisfies
\bea
S_M^{\pm}(\om) &=& 0, \qquad {\it if}\;\;\om\not\in\sg(M).  \lb{3.6}\\
S_{N_1(1,a)}^+(1) &=& \left\{
  \begin{array}{ll}1, &\quad {\rm if}\;\; a\ge 0, \cr
0, &\quad {\rm if}\;\; a< 0.  \end{array}\right.
  \lb{3.7}\eea

For any $M_i\in\Sp(2n_i)$ with $i=0$ and $1$, there holds
\be S^{\pm}_{M_0\diamond M_1}(\om) = S^{\pm}_{M_0}(\om) + S^{\pm}_{M_1}(\om),
    \qquad \forall\;\om\in\U. \lb{3.8}\ee\end{lemma}

For every $\ga\in\mathcal{P}_\tau(2n)\equiv\{\ga\in C([0,\tau],Sp(2n))\ |\ \ga(0)=I_{2n}\}$, we extend
$\ga(t)$ to $t\in [0,m\tau]$ for every $m\in\N$ by
\bea \ga^m(t)=\ga(t-j\tau)\ga(\tau)^j \qquad \forall\;j\tau\le t\le (j+1)\tau \;\;
               {\rm and}\;\;j=0, 1, \ldots, m-1, \lb{3.9}\eea
as in P.114 of \cite{lo1999}. As in \cite{LoZ} and \cite{lo2002}, we denote the Maslov-type indices of
$\ga^m$ by $(i(\ga,m),\nu(\ga,m))$.

Then the following decomposition theorem and precise iteration formula for symplectic paths are proved in \cite{lo1999}
and \cite{lo2000}.
\begin{theorem}\label{jqddgs}(cf. Theorem 1.8.10, Lemma 2.3.5 and Theorem 8.3.1 of \cite{lo2002})
For every $M\in\Sp(2n)$, there
exists a continuous path $f\in\Om^0(M)$ such that $f(0)=M$ and
\bea f(1)
&=& N_1(1,1)^{\dm p_-}\,\dm\,I_{2p_0}\,\dm\,N_1(1,-1)^{\dm p_+}
  \dm\,N_1(-1,1)^{\dm q_-}\,\dm\,(-I_{2q_0})\,\dm\,N_1(-1,-1)^{\dm q_+} \nn\\
&&\dm\,N_2(e^{\aa_{1}\sqrt{-1}},A_{1})\,\dm\,\cdots\,\dm\,N_2(e^{\aa_{r_{\ast}}\sqrt{-1}},A_{r_{\ast}})
  \dm\,N_2(e^{\bb_{1}\sqrt{-1}},B_{1})\,\dm\,\cdots\,\dm\,N_2(e^{\bb_{r_{0}}\sqrt{-1}},B_{r_{0}})\nn\\
&&\dm\,R(\th_1)\,\dm\,\cdots\,\dm\,R(\th_r)\dm\,D(\pm 2)^{\dm h},\lb{3.10}\eea
where
$\frac{\th_{j}}{2\pi}\in(0,1)$ for $1\le j\le r$; $N_2(e^{\aa_{j}\sqrt{-1}},A_{j})$'s
are nontrivial and $N_2(e^{\bb_{j}\sqrt{-1}},B_{j})$'s are trivial, and non-negative integers
$p_-, p_0, p_+,q_-, q_0, q_+,r,r_\ast,r_0,h$ satisfy
\be p_- + p_0 + p_+ + q_- + q_0 + q_+ + r + 2r_{\ast} + 2r_0 + h = n. \nn\ee

Let $\ga\in\P_{\tau}(2n)=\{\ga\in C([0,\tau],\Sp(2n))\,|\,\ga(0)=I\}$. Denote the basic normal form
decomposition of $M\equiv \ga(\tau)$ by (\ref{3.10}). Then we have
\bea i(\ga^m)
&=& m(i(\ga)+p_-+p_0-r ) + 2\sum_{j=1}^rE\left(\frac{m\th_j}{2\pi}\right) - r   \nn\\
&&  - p_- - p_0 - {{1+(-1)^m}\over 2}(q_0+q_+)
              + 2\sum_{j=1}^{r_{\ast}}\vf\left(\frac{m\aa_j}{2\pi}\right) - 2r_{\ast}. \nn\eea\end{theorem}

The common index jump theorem (cf. Theorem 4.3 of \cite{LoZ}) for symplectic paths established by Long
and Zhu in 2002 has become one of the main tools to study the multiplicity and stability problems of
closed solution orbits in Hamiltonian and symplectic dynamics. Recently, the following enhanced version
of it has been obtained by Duan, Long and Wang in \cite{DLW2}, which will play an important role in the
proofs in Section 6.
\begin{theorem}\label{ECIJT}(cf. Theorem 3.5 of \cite{DLW2}) ({\bf The enhanced common index jump theorem for
symplectic paths})  Let $\gamma_k\in\mathcal{P}_{\tau_k}(2n)$ for $k=1,\cdots,q$ be a finite
collection of symplectic paths. Let $M_k=\ga(\tau_k)$. We extend $\ga_k$ to $[0,+\infty)$ by (\ref{3.9})
inductively. Suppose
\be  \hat{i}(\ga_k,1) > 0, \qquad \forall\ k=1,\cdots,q.  \lb{3.11}\ee
Then for every integer $\bar{m}\in \N$, there exist infinitely many $(q+1)$-tuples
$(N, m_1,\cdots,m_q) \in \N^{q+1}$ such that for all $1\le k\le q$ and $1\le m\le \bar{m}$, there holds
\bea
\nu(\ga_k,2m_k-m) &=& \nu(\ga_k,2m_k+m) = \nu(\ga_k, m),   \lb{3.12}\\
i(\ga_k,2m_k+m) &=& 2N+i(\ga_k,m),                         \lb{3.13}\\
i(\ga_k,2m_k-m) &=& 2N-i(\ga_k,m)-2(S^+_{M_k}(1)+Q_k(m)),  \lb{3.14}\\
i(\ga_k, 2m_k)&=& 2N -(S^+_{M_k}(1)+C(M_k)-2\Delta_k),     \lb{3.15}\eea
where
\be \Delta_k = \sum_{0<\{m_k\th/\pi\}<\delta}S^-_{M_k}(e^{\sqrt{-1}\th}),\qquad
 Q_k(m) = \sum_{e^{\sqrt{-1}\th}\in\sg(M_k),\atop \{\frac{m_k\th}{\pi}\}
                   = \{\frac{m\th}{2\pi}\}=0}S^-_{M_k}(e^{\sqrt{-1}\th}). \lb{3.16}\ee
More precisely, by (4.10), (4.40) and (4.41) in \cite{LoZ} , we have
\bea m_k=\left(\left[\frac{N}{\bar{M}\hat i(\gamma_k, 1)}\right]+\chi_k\right)\bar{M},\quad 1\le k\le q,\lb{3.17}\eea
where $\chi_k=0$ or $1$ for $1\le k\le q$ and $\frac{\bar{M}\theta}{\pi}\in\Z$
whenever $e^{\sqrt{-1}\theta}\in\sigma(M_k)$ and $\frac{\theta}{\pi}\in\Q$
for some $1\le k\le q$.  Furthermore, given $M_0\in\N$,
by (iv) of Remark 3.6 of \cite{DLW2}, we may
further require  $M_0|N$, and by (4.20) in Theorem 4.1 of \cite{LoZ},
for any $\epsilon>0$, we can choose $N$ and $\{\chi_k\}_{1\le k\le q}$ such that
\bea \left|\left\{\frac{N}{\bar{M}\hat i(\gamma_k, 1)}\right\}-\chi_k\right|<\epsilon,\quad 1\le k\le q.\lb{3.18}\eea
In particular, if $\nu(\ga_k, m)=0$ for all $1\leq k\leq q$ and $m\in\N$, by (\ref{3.6}) and (\ref{3.16}) we have $S^+_{M_k}(1)=0$ and $Q_k(m)=0$, $\forall\ m\ge 1$,
then (\ref{3.13})-(\ref{3.15}) have the following simple forms
\bea
i(\ga_k,2m_k+m) &=& 2N+i(\ga_k,m),                         \lb{3.19}\\
i(\ga_k,2m_k-m) &=& 2N-i(\ga_k,m),  \lb{3.20}\\
i(\ga_k, 2m_k)&=& 2N -C(M_k)+2\Delta_k.     \lb{3.21}\eea
\end{theorem}

\setcounter{equation}{0}
\section{Katok's metrics on spheres and compact space forms}
In this section, we compute the $S^1$-equivariant Betti number sequence of
$(\Lambda_e M, \Lambda^0 M)$ via Katok's famous metrics on $S^{n}$,
where $M=S^{2n+1}/ \Gamma$, $\Lambda^0 M  =\{{\rm constant\;point\;curves\;in\;}M\}\cong M$
and $\Lambda_e M$ is the contractible component of the free loop space on $M$.

In 1973, Katok \cite{Katok1973} constructed his famous irreversible Finsler metrics on $S^n$
which possess only finitely many distinct prime closed geodesics.
His examples were further studied closely by Ziller \cite{Ziller1982} in 1982,
from which we borrow most of the notations.

For $S^{2n-1}$ with the canonical Riemannian metric $g$, any closed one-parameter group of
isometries is conjugate to a diagonal
matrix \bea \phi_t=\diag(R(pt/p_1),\cdots, R(pt/p_{n})),\nn\eea
where $p_i\in\Z, p=p_1\cdots p_{n}$, the $p_i$ are relatively prime and $R(\omega)$ is
a rotation in $\R^2$ with angle $\omega$.
For $S^{2n}$ with the canonical Riemannian metric $g$, the
same is true if the matrix is enlarged by one row and one column with a 1 in the diagonal.
Let $TS^n$ and $T^*S^n$ be its tangent bundle and cotangent bundle respectively.
Define $H_0, H_1: T^*S^n\rightarrow \R$ by
\bea H_0(x)=\|x\|_*\quad and \quad H_1(x)=x(V), \forall x\in T^*S^n,\nn\eea
where $\|\cdot\|_*$ denotes the dual norm of $g$ and $V$ is the vector field generated by $\phi_t$. Let
\bea H_\alpha= H_0 +\alpha H_1 \quad for \quad \alpha\in (0, 1).\nn\eea
Then $\frac{1}{2}H^2_\alpha$ is homogeneous of degree two and the Legendre transform
\bea L_{\frac{1}{2}H^2_\alpha}=D_F\left(\frac{1}{2}H^2_\alpha\right): T^*S^n\rightarrow  TS^n\nn\eea
is a global diffeomorphism. Hence,\bea N_\alpha=H_\alpha\circ L_{\frac{1}{2}H^2_\alpha}^{-1}\nn\eea
defines a Finsler metric on $S^{n}$. Since $H_\alpha(-x) \neq H_\alpha(x)$, $N_\alpha$ is not reversible. It was
proved that $(S^n, N_\alpha)$ with $\alpha\in(0, 1) \setminus \Q$ possesses precisely $2[\frac{n +1}{2}]$ distinct
prime closed geodesics and all of them are irrationally
elliptic(cf. Katok \cite{Katok1973} and pp.137-139 of Ziller \cite{Ziller1982} for more details).

Consider a finite group $\Gamma$ which acts freely and isometrically on $(S^{n}, g)$.
For $(S^{2n-1}, g)$, any element $h$ of $\Gamma$
has the form \bea h=\diag(R(\alpha_1),\cdots, R(\alpha_{n}))£¬\nn\eea
for some $\alpha_i\in[0,2\pi)$, where $1\leq i\leq n$.
Note that the only non-trivial group which acts freely on $S^{2n}$ is $\Z_2$(cf. P.5 of \cite{Tai2016}),
then $h$ is the identity or the antipodal map for $(S^{2n}, g)$. Thus $h\circ\phi_t=\phi_t\circ h$ for any $(S^n, g)$.
\begin{lemma}\label{iso} For any element $h\in\Gamma$, $h: (S^n, N_\alpha)\rightarrow
(S^n, N_\alpha)$ is an isometry.\end{lemma}
\Proof For any $x\in T_{p}^*S^n$, we have $H_\alpha\circ (h^{-1})^*(x)=H_\alpha(x)$.
In fact, we have $(h^{-1})^*(x)\in T_{h(p)}^*S^n$ and
\bea H_\alpha\circ (h^{-1})^*(x)&=&\|(h^{-1})^*(x)\|_*+\alpha ((h^{-1})^*(x))(V_{h(p)})\nn\\
&=&\|x\|_*+\alpha x((h^{-1})_*(V_{h(p)}))\nn\\
&=&\|x\|_*+\alpha x(V_{p})\nn\\&=&H_\alpha(x)
\lb{hbb}\eea
where the second identity is due to the fact that $h$ is isomorphic for the canonical metric $\|\cdot\|_*$
and the third identity is due to the fact that $h_*(V_p)=V_{h(p)}$ which follows by the definition of $V$
and $h\circ\phi_t=\phi_t\circ h$.

To prove that $h$ is an isometry of $(S^n, N_\alpha)$, i.e.,
\bea  H_\alpha \circ L_{\frac{1}{2}H^2_\alpha}^{-1}\circ h_*(X)=
H_\alpha \circ L_{\frac{1}{2}H^2_\alpha}^{-1}(X),\forall X\in T_{p}S^n,\nn\eea
it is sufficient by (\ref{hbb}) to prove
\bea L_{\frac{1}{2}H^2_\alpha}^{-1}\circ h_*(X)=
(h^{-1})^* \circ L_{\frac{1}{2}H^2_\alpha}^{-1}(X),\forall X\in T_{p}S^n.\lb{jh}\eea
In fact, (\ref{jh}) is equivalent to
\bea h_*\circ L_{\frac{1}{2}H^2_\alpha}(x)=
L_{\frac{1}{2}H^2_\alpha}\circ (h^{-1})^*(x),\forall x\in T^*_{p}S^n.\lb{jhbx}\eea
Note that for any $x\in T^*_{p}S^n$,
\bea L_{\frac{1}{2}H^2_\alpha}(x)&=&D_F\left(\frac{1}{2}H^2_\alpha\right)(x)\nn\\&=&H_\alpha(x)\cdot D_F(H_\alpha)(x)
\nn\\&=&H_\alpha(x)(*x/\|*x\|+\alpha V_{p})
\nn\eea
where $*x$ is the canonical identification of $x$(cf. p. 143 in \cite{Ziller1982}). Thus we get
\bea  L_{\frac{1}{2}H^2_\alpha}\circ (h^{-1})^*(x)&=&H_\alpha\circ
(h^{-1})^*(x)\cdot(*((h^{-1})^*x)/\|*((h^{-1})^*x)\|+\alpha V_{h(p)})\nn\\
&=&H_\alpha(x)(*((h^{-1})^*x)/\|*((h^{-1})^*x)\|+\alpha V_{h(p)})\nn\\&=&H_\alpha(x)
(h_*(*x)/\|h_*(*x)\|+\alpha V_{h(p)})\nn\\&=&H_\alpha(x)
(h_*(*x)/\|*x\|+\alpha V_{h(p)})\nn\\&=&h_*\circ L_{\frac{1}{2}H^2_\alpha}(x),\nn\eea
where the second identity is due to (\ref{hbb}), the third identity is due to the fact that $*((h^{-1})^*x)=h_*(*x)$
, the fourth identity is due to the fact that $h$ is isomorphic for the canonical metric $\|\cdot\|$ and
the fifth identity is due to the fact that $h_*(V_p)=V_{h(p)}$.
Then (\ref{jhbx}), and then (\ref{jh}) as well as Lemma \ref{iso} are proved.\hfill$\Box$

By Lemma \ref{iso}, we can endow $S^n/ \Gamma$ a Finsler metric induced by $N_\alpha$,
which is still denoted by $N_\alpha$ for simplicity.
Therefore the natural projection
\bea \pi:(S^n, N_\alpha)\rightarrow (S^n/ \Gamma, N_\alpha)\nn\eea
is locally isometric.

In order to compute the $S^1$-equivariant Betti number sequence of
$(\Lambda_e M, \Lambda^0 M)$ for $M=S^{2n+1}/ \Gamma$, we review the $S^1$-equivariant Betti
numbers of the free loop space pair $(\Lambda S^{2n+1}, \Lambda^0 S^{2n+1})$.
\begin{lemma}\label{qmloop} (cf. Theorem 2.4 and Remark 2.5 of \cite{Rad1989} and \cite{Hin1984}, cf.
also Lemma 2.5 of \cite{DL2010})  Let $(S^{2n+1},F)$ be a
$(2n+1)$-dimensional Finsler sphere. Then the $S^1$-equivariant Betti numbers
of $(\Lambda S^{2n+1}, \Lambda^0 S^{2n+1})$ are given by

\bea b_j\equiv\rank H_j(\Lm S^{2n+1}/S^1,\Lm^0 S^{2n+1}/S^1;\Q)  =\left\{\begin{array}{ll}
    2,&\quad {\it if}\quad j\in \K\equiv \{2nk\,|\,2\le k\in\N\},  \\
    1,&\quad {\it if}\quad j\in \{2n+2k\,|\,k\in\N_0\}\bs\K,  \\
    0 &\quad {\it otherwise}. \\
\end{array}
\right. \lb{qmbds}\eea
\end{lemma}

Combining Lemma \ref{qmloop} with the fact that $\pi:(S^n, N_\alpha)\rightarrow (S^n/ \Gamma, N_\alpha)$
is locally isometric, we can give the $S^1$-equivariant Betti numbers of
$(\Lambda_e M, \Lambda^0 M)$ as follows.
\begin{proposition}\label{ksfztd}   For $M=S^{2n+1}/ \Gamma$, the $S^1$-equivariant Betti numbers
of $(\Lambda_e M, \Lambda^0 M)$ are given by
\bea \beta_j\equiv\rank H_j(\Lm_e M/S^1,\Lambda^0 M/S^1;\Q)  =\left\{\begin{array}{ll}
    2,&\quad {\it if}\quad j\in \K\equiv \{2nk\,|\,2\le k\in\N\},  \\
    1,&\quad {\it if}\quad j\in \{2n+2k\,|\,k\in\N_0\}\bs\K,  \\
    0 &\quad {\it otherwise}. \\
\end{array}
\right. \lb{ksfzbds}\eea
and the average $S^1$-equivariant Betti number of $(\Lambda_e M, \Lambda^0 M)$ satisfies
\be  \bar{B}(\Lambda_e M, \Lambda^0 M;\Q)\equiv\lim_{q\to+\infty}\frac{1}{q}
\sum_{k=0}^{q}(-1)^k\beta_k= \frac{n+1}{2n}.   \lb{aB.1}\ee
\end{proposition}

\Proof Since the Betti numbers $\beta_j$ are topological invariants of $M$, they are independent
of the choice of the Finsler metric $F$ on it. To estimate them, it suffices to choose a special Finsler metric
$F=N_\alpha$ for $\alpha\in(0, 1) \setminus \Q$, i.e., the Katok metrics.

Since all prime closed geodesics on $(S^{2n+1} ,N_\alpha)$ are irrationally elliptic,
denoted by $\{c_j\}_{j=1}^{2n+2}$.
By the decomposition in Theorem \ref{jqddgs}, the linearized Poincar\'{e} map $P_{c_j}$
can be connected to $f_{c_j}(1)$ in $\Om^0(P_{c_j})$ satisfying
\bea f_{c_j}(1)&=&R(\th_{j1})\,\dm\,\cdots\,\dm\,R(\th_{j(2n)}),
\quad 1\le j\le 2n+2, \lb{wltyfj}\eea
where $\frac{\th_{jk}}{2\pi}\in\Q^c$ for $1\le j\le 2n+2$ and $1\le k\le 2n$.
Then by Theorem \ref{jqddgs}, we obtain
\bea i(c_j)&\in& 2\Z,\qquad 1\le j\le 2n+2,\lb{zbjex}\\
i(c_j^m)&=& m(i(c_j)-2n) + 2\sum_{k=1}^{2n}\left[\frac{m\th_{jk}}{2\pi}\right]
+2n,\quad \nu(c_j^m)=0,\qquad\forall\ m\ge 1, \lb{zbgs}\eea
where (\ref{zbjex}) follows from Theorem 8.1.7 of \cite{lo2002} and the symplectic additivity of iterated indices.

By (\ref{zbjex}) and (\ref{zbgs}), there holds
\bea i(c_j^m)\in 2\Z,\qquad \forall\ 1\le j\le 2n+2,\quad m\ge 1. \lb{jexxt}\eea

Define
\bea M_p=\sum_{j=1}^{2n+2} M_p(j)\equiv\sum_{j=1}^{2n+2}\#\{m\ge 1\ |\ i(c_j^m)
=p, \;\ol{C}_p(E, c_j^m)\not= 0\},\quad p\in\Z.\nn\eea
Then the following Morse inequality (cf. Theorem I.4.3 of \cite{Cha}) holds
\bea M_p&\ge& b_p,\lb{mes1}\\
M_p - M_{p-1} + \cdots +(-1)^{p}M_0
&\ge& b_p - b_{p-1}+ \cdots + (-1)^{p}b_0, \quad\forall\ p\in\N_0.\lb{mes2}\eea
where the $b_j's$ are given in Lemma \ref{qmloop}.

By Lemma \ref{Rad1992} and (\ref{jexxt}) we obtain
\bea M_p=\sum_{j=1}^{2n+2} M_p(j)=
\sum_{j=1}^{2n+2}\#\{m\ge 1\ |\ i(c_j^m)=p\},\qquad \forall\ p\in\Z, \lb{mess}\eea
which together with (\ref{jexxt}), yields
\be M_{p}=b_{p}=0, \qquad \forall\ p=1\ (\mod 2),\ p\in\N_0.{\lb{mbxd}}\ee

Then by the Morse inequality (\ref{mes2}), we obtain
\be M_{p}=b_{p}, \qquad \forall\ p=0\ (\mod 2),\ p\in\N_0.{\lb{mbxd2}}\ee

On the other hand, the set of contractible closed geodesics on
$(M, N_\alpha)=(S^{2n+1}/ \Gamma, N_\alpha)$ is
$\{c_j^m\mid 1\le j\le 2n+2, m\in\N\}$ since the closed geodesics $c_j's$ found
in \cite{Ziller1982} are great circles of $S^{2n+1}$.
Then $c_j$ is a contractible minimal closed geodesic on $(S^{2n+1}/ \Gamma, N_\alpha)$.
Note that (\ref{jexxt}) yields
\bea i(c_j^m)-i(c_j)\in 2\Z,\qquad \forall\ 1\le j\le 2n+2,\quad m\ge 1. \lb{jexxt2}\eea
Combining (\ref{jexxt2}) with Proposition \ref{Rad1992'}, we know the Morse type numbers of the contractible component
of the free loop space on $(S^{2n+1}/ \Gamma, N_\alpha)$ are just the Morse type numbers $M_p$'.
Then we can apply Morse inequality to the contractible component
of the free loop space on $(S^{2n+1}/ \Gamma, N_\alpha)$
and obtain \bea M_p&\ge& \beta_p,\lb{ksmes1}\\
M_p - M_{p-1} + \cdots +(-1)^{p}M_0
&\ge& \beta_p - \beta_{p-1}+ \cdots + (-1)^{p}\beta_0, \quad\forall\ p\in\N_0.\lb{ksmes2}\eea
Combining (\ref{mbxd}) with (\ref{ksmes1}), we have
\be \beta_{p}=0, \qquad \forall\ p=1\ (\mod 2),\ p\in\N_0.{\lb{lmbxd}}\ee

Then by the Morse inequality (\ref{ksmes2}), we obtain
\be M_{p}=\beta_{p}, \qquad \forall\ p=0\ (\mod 2),\ p\in\N_0, {\lb{lmbxd2}}\ee
which together with (\ref{mbxd2}) yields
\be \beta_{p}=b_{p}, \qquad \forall\ p=0\ (\mod 2),\ p\in\N_0.{\lb{bdsjg}}\ee
It together with (\ref{lmbxd}) and Lemma \ref{qmloop} completes our proofs. \hfill$\Box$

\setcounter{equation}{0}
\section{Resonance identity for contractible closed geodesics on $(S^{2n+1}/ \Gamma,F)$}

In this section, we apply Proposition \ref{ksfztd} to obtain the resonance identity for homologically
visible contractible minimal closed geodesics on a Finsler $M=(S^{2n+1}/ \Gamma,F)$ claimed in Theorem \ref{Thm1.1},
provided the number of all distinct contractible minimal closed geodesics on $M$ is finite.
Note that if the number of all distinct contractible minimal closed geodesics on $M$ is finite,
then the total number of all distinct closed geodesics on $M$ is finite.

Firstly we have
\begin{definition}\label{def-hv} Let $(M,F)$ be a compact Finsler manifold. A contractible closed geodesic $c$ on $M$
is homologically visible, if there exists an integer $k\in\Z$ such that $\bar{C}_k(E,c; [e]) \not= 0$. We denote
by $\CG_{\hv}(M,F)$ the set of all distinct homologically visible contractible minimal closed geodesics on $(M,F)$.
\end{definition}

{\bf Proof of Theorem \ref{Thm1.1}.} Recall that we denote the homologically visible contractible
minimal closed geodesics by $\CG_{\hv}(M)=\{c_1, \ldots, c_r\}$ for some integer $r>0$
when the number of distinct contractible minimal
closed geodesics on $M=(S^{2n+1}/ \Gamma,F)$ is finite.

Note also that
we have $\hat{i}(c_j)>0$ for all
$1\le j\le r$, otherwise there is a homologically visible contractible closed geodesic
$c$ on $M$ satisfying $\hat{i}(c)=0$. Then $i(c^m) =0$ for all $m\in \N$
by Bott iteration formula and $c$ must be an absolute minimum
of the energy functional $E$ in its free homotopy class,
since otherwise there would exist infinitely many distinct closed geodesics on $M$ by Theorem 3
on p.385 of \cite{BK1983}. It also follows that this homotopy class must be non-trivial,
which contradicts to that $c$ is contractible.

Denote the contractible component of the free loop space on $M$ by $\Lambda_e M$.
Then the set of contractible closed geodesics on $M=(S^{2n+1}/ \Gamma, F)$ is
$\{c_j^m\mid 1\le j\le r, m\in\N\}$.

Let \bea M_q\equiv M_q(\Lambda_e M, \Lambda^0 M) =
\sum_{1\le j\le r,\; m\ge 1}\dim{\ol{C}}_q(E, c^m_j; [e]),\quad q\in\Z.\lb{dfmess}\eea
The Morse series of $(\Lambda_e M, \Lambda^0 M)$ is defined by
\be  M(t) = \sum_{h=0}^{+\infty}M_ht^h.  \label{wh}\ee

{\bf Claim 1.} {\it $\{m_h\}$ is a bounded sequence.}

In fact, by (\ref{CGap2}), we have
\bea M_h=\sum_{j=1}^{r}\sum_{m=1}^{n_j}\sum_{l=0}^{4n}k^{\epsilon(c_j^{m})}_{l}(c_{j}^{m})
            \;{}^{\#}\left\{s\in\mathbb{N}\cup\{0\}\mid h-i(c_{j}^{m+sn_j})=l\right\}, \label{wh2}\eea
and by Theorems 10.1.2 of \cite{lo2002},
we have $|i(c_{j}^{m+sn_j})-(m+sn_j)\hat{i}(c_{j})|\le 2n$, then
\bea
&&{}^{\#}\left\{s\in\mathbb{N}\cup\{0\}\mid h-i(c_{j}^{m+sn_j})=l\right\} \nn\\
&&\qquad = \;{}^{\#}\left\{s\in\mathbb{N}\cup\{0\}\mid l+i(c_{j}^{m+sn_j})=h,\;
           |i(c_{j}^{m+sn_j})-(m+sn_j)\hat{i}(c_{j})|\le 2n\right\}  \nn\\
&&\qquad \le \;{}^{\#}\left\{s\in\mathbb{N}\cup\{0\}\mid 2n\geq |h-l-(m+sn_{j})\hat{i}(c_{j})|  \right\}  \nn\\
&&\qquad = \;{}^{\#}\left\{s\in\mathbb{N}\cup\{0\}\mid \frac{h-l-2n-m\hat{i}(c_j)}{n_j\hat{i}(c_j)} \leq s
       \le \frac{h-l+2n-m\hat{i}(c_j)}{n_j\hat{i}(c_j)}\right\} \nn\\
&&\qquad \le \frac{4n}{n_j\hat{i}(c_j)} + 1.  \label{wh3}\eea
Hence Claim 1 follows by (\ref{wh2}) and (\ref{wh3}).

We now use the method
in the proof of Theorem 5.4 of \cite{LW2007} to estimate
$$   M^{q}(-1) = \sum_{h=0}^{q}M_h(-1)^h.   $$

By (\ref{wh}) and (\ref{CGap2}) we obtain
\bea M^{q}(-1)
&=& \sum_{h=0}^{q}M_{h}(-1)^{h}  \nn\\
&=& \sum_{j=1}^{r}\sum_{m=1}^{n_j}\sum_{l=0}^{4n}\sum_{h=0}^{q}(-1)^{h}k^{\epsilon(c_j^{m})}_{l}(c_{j}^{m})
            \;{}^{\#}\left\{s\in\mathbb{N}\cup\{0\}\mid h-i(c_{j}^{m+sn_j})=l\right\}  \nn\\
&=& \sum_{j=1}^{r}\sum_{m=1}^{n_j}\sum_{l=0}^{4n}(-1)^{l+i(c^m_{j})}k^{\epsilon(c_j^{m})}_{l}(c_{j}^{m})
            \;{}^{\#}\left\{s\in\mathbb{N}\cup\{0\}\mid l+i(c_{j}^{m+sn_j})\le q\right\}.  \nn\eea
On the one hand, we have
\bea
&&{}^{\#}\left\{s\in\mathbb{N}\cup\{0\}\mid l+i(c_{j}^{m+sn_j})\le q\right\} \nn\\
&&\qquad = \;{}^{\#}\left\{s\in\mathbb{N}\cup\{0\}\mid l+i(c_{j}^{m+sn_j})\le q,\;
           |i(c_{j}^{m+sn_j})-(m+sn_j)\hat{i}(c_{j})|\le 2n\right\}  \nn\\
&&\qquad \le \;{}^{\#}\left\{s\in\mathbb{N}\cup\{0\}\mid 0\le (m+sn_{j})\hat{i}(c_{j})\le q-l+2n  \right\}  \nn\\
&&\qquad = \;{}^{\#}\left\{s\in\mathbb{N}\cup\{0\}\mid 0 \leq s
       \le \frac{q-l+2n-m\hat{i}(c_j)}{n_j\hat{i}(c_j)}\right\} \nn\\
&&\qquad \le \frac{q-l+2n}{n_j\hat{i}(c_j)} + 1. \nn\eea
On the other hand, we have
\bea
&&{}^{\#}\left\{s\in\mathbb{N}\cup\{0\}\mid l+i(c_{j}^{m+sn_j})\le q \right\} \nn\\
&&\qquad = \;{}^{\#}\left\{s\in\mathbb{N}\cup\{0\}\mid l+i(c_{j}^{m+sn_j})\le q,\;
           |i(c_{j}^{m+sn_j})-(m+sn_{j})\hat{i}(c_{j})|\le 2n\right\}  \nn\\
&&\qquad \ge \;{}^{\#}\left\{s\in\mathbb{N}\cup\{0\}\mid i(c_{j}^{m+sn_j})
       \le (m+sn_{j})\hat{i}(c_{j})+2n\le q-l \right\} \nn\\
&&\qquad \ge \;{}^{\#}\left\{s\in\mathbb{N}\cup\{0\}\mid 0 \le s
       \le \frac{q-l-2n-m\hat{i}(c_{j})}{n_j\hat{i}(c_{j})} \right\}  \nn\\
&&\qquad \ge \frac{q-l-2n}{n_j\hat{i}(c_{j})} - 1.  \nn\eea
Thus we obtain
$$  \lim_{q\to+\infty}\frac{1}{q}M^{q}(-1)
  = \sum_{j=1}^{r}\sum_{m=1}^{n_j}\sum_{l=0}^{4n}(-1)^{l+i(c_{j}^m)}
  k^{\epsilon(c_j^{m})}_{l}(c_{j}^m)\frac{1}{n_{j}\hat{i}(c_{j})}
              = \sum_{j=1}^{r}\frac{\hat{\chi}(c_j)}{\hat{i}(c_j)}.  $$
Since $m_{h}$ is bounded, by Morse inequality (cf. Theorem I.4.3 of \cite{Cha}) we then obtain
$$  \lim_{q\to+\infty}\frac{1}{q}M^{q}(-1)
          = \lim_{q\to+\infty}\frac{1}{q}\sum_{k=0}^{q}(-1)^k\beta_k = \bar{B}(\Lambda_e M, \Lambda^0 M;\Q),  $$
Thus by (\ref{aB.1}) we get
$$  \sum_{j=1}^{r}\frac{\hat{\chi}(c_j)}{\hat{i}(c_j)} =
    \frac{n+1}{2n}, $$
which proves (\ref{reident1}) of Theorem \ref{Thm1.1}.
For the special case when each $c_{j}^{m}$
is non-degenerate with $1\leq j\leq r$ and $m\in\mathbb{N}$,
we have $n_{j}=2$ and $k_{l}^{\epsilon(c_j^{m})}(c_{j}^{m})=1$ when $l=0$ and $i(c_j^{m})-i(c_j)\in 2\Z$,
and $k_{l}^{\epsilon(c_j^{m})}(c_{j}^{m})=0$ for all other
$l\in\Z$. Then (\ref{reident1}) has the following simple form
\be  \sum_{j=1}^{r}\left(\frac{(-1)^{i(c_{j})}k_0^{\epsilon(c_j)}(c_{j})+(-1)^{i(c_{j}^{2})}
k_0^{\epsilon(c_j^{2})}(c_{j}^{2})}{2}\right)\frac{1}{\hat{i}(c_{j})}=
    \frac{n+1}{2n},  \nn\ee
which proves (\ref{breident1}) of Theorem \ref{Thm1.1}. $\hfill\Box$

\setcounter{equation}{0}
\section{Proof of Theorem 1.2}

In order to prove Theorem 1.2, let $M=(S^{2n+1}/ \Gamma,F)$ with a bumpy, irreversible Finsler metric F.
We make the following assumption

{\bf (FCCG)} {\it Suppose that there exist only finitely many contractible minimal closed geodesics $\{c_k\}_{k=1}^q$
with $i(c_k)\geq 2$ for $k=1,2,\cdots,q$ on $M$.}

Since the flag curvature $K$ of $(S^{2n+1}/ \Gamma, F)$ satisfies
$\left(\frac{\lambda}{\lambda+1}\right)^2<K\le 1$ by assumption,
then every contractible closed geodesic $c_k$ must satisfy
\bea i(c_k)\ge 2n\geq 2, \nn\eea
by Theorem 3 and Lemma 3 of \cite{Rad2004}. Thus (FCCG) holds.

{\bf Step 1} {\it There exist at least $2n$ distinct non-hyperbolic contractible closed geodesics, all of which
possess even Morse indices.}

Since by the assumption (FCCG), there exist only finitely many contractible minimal closed geodesics
on the bumpy manifold $M$, any closed geodesic $c_k$ among $\{c_k \}^q_{k=1}$ must have
positive mean index (cf. Section 5), i.e.,
\bea \hat{i}(c_k)>0,\quad 1\leq k\leq q,\lb{pmind}\eea
which implies that $i(c_k^m)\rightarrow +\infty$ as $m\rightarrow +\infty$.
So the positive integer $\bar{m}$ defined by
\bea  \bar{m}=\max_{1\le k\le q}\left\{\min\{m_0\in\N\ |\ i(c_k^{m+l})\ge i(c_k^l), \quad
    \forall\ l\ge 1, m\geq m_0\}\right\}  \lb{mba}\eea
is well-defined and finite.

For the integer $\bar{m}$ defined in (\ref{mba}), it follows from (\ref{3.19})-(\ref{3.21}) of Theorem \ref{ECIJT}
that there exist infinitely many $q+1$-tuples $(N, m_1, \cdots, m_q)\in\N^{q+1}$ such that for any
$1\le k\le q$, there holds
\bea
\bar{m}+2 &\le& \min\{2m_k,\ 1\le k\le q\},\lb{ECIJT1}\\
i(c_k^{2m_k-m}) &=& 2N-i(c_k^m),\quad 1\le m\le\bar{m}, \lb{ECIJT2}\\
i(c_k^{2m_k}) &=& 2N-C(M_k)+2\Delta_k,\lb{ECIJT3}\\
i(c_k^{2m_k+m})&=& 2N+i(c_k^m),\quad 1\le m\le\bar{m},\lb{ECIJT4}\eea
where $M_k=P_{c_k}\in Sp(4n)$ is the linearized Poincar$\acute{e}$ map of the contractible minimal closed geodesic
$c_k$. Here note that in the bumpy case, we have $\nu(c_k^m)=0$ for all $1\leq k\leq q$ and $m\in\N$.

On one hand, there holds $i(c^m_k)\geq i(c_k)$ for any $m\geq 1$ by the Bott-type formulae (cf.
\cite{Bott1956} and Theorem 9.2.1 of \cite{lo2002}). Thus by (\ref{ECIJT2}), (\ref{ECIJT4})
and the assumption $i(c_k)\geq 2$ for
$1 \leq k \leq q$ in the theorem, it yields
\bea
i(c_k^{2m_k-m}) &=& 2N-i(c_k^m)\leq 2N-i(c_k)\leq 2N-2,\quad 1\le m\le\bar{m}, \lb{mbzbgj1}\\
i(c_k^{2m_k+m})&=& 2N+i(c_k^m)\geq 2N+i(c_k)\geq 2N+2,\quad 1\le m\le\bar{m},\lb{mbzbgj1}\eea
On the other hand, by the definition (\ref{mba}) of $\bar{m}$ and (\ref{ECIJT2}), (\ref{ECIJT4}),
for $2m_k>m\ge \bar{m}+1$, we obtain
\bea
i(c_k^{2m_k-m}) &\leq& i(c_k^{2m_k-1})= 2N-i(c_k)\leq 2N-2,\lb{mbzbgj3}\\
i(c_k^{2m_k+m})&\geq& i(c_k^{2m_k+1})= 2N+i(c_k)\geq 2N+2.\lb{mbzbgj4}\eea
In summary, by (\ref{ECIJT2})-(\ref{mbzbgj4}), for $1\leq k\leq q$, we have proved
\bea
i(c_k^{2m_k-m}) &\leq& 2N-2,\quad 1\leq m< 2m_k, \lb{zbgj1}\\
i(c_k^{2m_k}) &=& 2N-C(M_k)+2\Delta_k,\lb{zbgj2}\\
i(c_k^{2m_k+m})&\geq& 2N+2,\quad \forall m\geq 1.\lb{zbgj3}\eea

{\bf Claim 1:} {\it For $N\in\N$ in Theorem \ref{ECIJT} satisfying (\ref{zbgj1})-(\ref{zbgj3})
and $2N\bar{B}(\Lambda_e M, \Lambda^0 M;\Q)\in\N$, we have
\bea \sum_{1\le k\le q} 2m_k\hat{\chi}(c_k)=2N\bar{B}(\Lambda_e M, \Lambda^0 M;\Q).\lb{claim1}\eea}

In fact, let $\ep<\frac{1}{1+2\bar{M}\sum_{1\le k\le q}|\hat{\chi}(c_k)|}$, by Theorem \ref{Thm1.1} and
(\ref{3.17})-(\ref{3.18}), it yields
\bea \left|2N\bar{B}(\Lambda_e M, \Lambda^0 M;\Q)-\sum_{k=1}^q 2m_k\hat{\chi}(c_k)\right|
&=& \left|\sum_{k=1}^q\frac{2N\hat{\chi}(c_k)}{\hat{i}(c_k)}-\sum_{k=1}^q 2\hat{\chi}(c_k)
    \left(\left[\frac{N}{\bar{M}\hat{i}(c_k)}\right]+\chi_k\right)\bar{M}\right|\nn\\
&\le& 2\bar{M}\sum_{k=1}^q |\hat{\chi}(c_k)|\left|\left\{\frac{N}{\bar{M}\hat{i}(c_k)}\right\}-\chi_k\right|.\nn\\
&<& 2\bar{M}\ep\sum_{k=1}^q|\hat{\chi}(c_k)|\nn\\
&<& 1,\nn\eea
which proves Claim 1 since each $2m_k\hat{\chi}(c_k)$ is an integer.

\medskip

Now by Proposition \ref{Rad1992'}, it yields
\bea
&& \sum_{m=1}^{2m_k} (-1)^{i(c_k^m)}\dim \overline{C}_{i(c_k^{m})}(E,c_k^m; [e])\nn\\
&& = \sum_{i=0}^{m_k-1} \sum_{m=2i+1}^{2i+2} (-1)^{i(c_k^m)} \dim
\overline{C}_{i(c_k^{m})}(E,c_k^m; [e])\nn\\
&& = \sum_{i=0}^{m_k-1} \sum_{m=1}^{2} (-1)^{i(c_k^m)} \dim \overline{C}_{i(c_k^{m})}(E,c_k^m; [e])\nn\\
&& = m_k \sum_{m=1}^{2} (-1)^{i(c_k^m)} \dim \overline{C}_{i(c_k^{m})}(E,c_k^m; [e])\nn\\
&& = 2m_k\hat{\chi}(c_k),\qquad \forall\ 1\le k\le q, \lb{pjzbdsbx}\eea
where the second equality follows from Proposition \ref{Rad1992'} and the fact $i(c^{m+2}_k)-i(c^m_k)\in 2\Z$ for
all $m\geq 1$ from Theorem 3.1, and the last equality follows from
Proposition \ref{Rad1992'}, (\ref{breident1}) and (\ref{CGht1}).

So, for $1\le k\le q$, by (\ref{zbgj1}), (\ref{zbgj3}) and Proposition \ref{Rad1992'}, we know that all $c^{2m_k-m}_k$'s
with $2m_k>m \geq 1$ only have contributions to the alternating sum
$\sum^{2N-2}_{p=0}(-1)^pM_p$, and all $c^{2m_k+m}_k$'s with $m \geq 1$
have no contributions to $\sum^{2N+1}_{p=0}(-1)^pM_p$.

Thus for the Morse-type numbers $M_p$'s in (\ref{dfmess}), by (\ref{pjzbdsbx}) we have
\bea \sum_{p=0}^{2N+1}(-1)^p M_p
&=& \sum_{k=1}^{q}\ \sum_{1\le m\le 2m_k \atop i(c_k^{m})\le 2N+1}
    (-1)^{i(c_k^m)}\dim \overline{C}_{i(c_k^{m})}(E,c_k^m; [e])\nn\\
&=& \sum_{k=1}^{q}\ \sum_{m=1}^{2m_k} (-1)^{i(c_k^m)}\dim \overline{C}_{i(c_k^{m})}(E,c_k^m; [e])\nn\\
& & -\sum_{1\le k\le q \atop i(c_k^{2m_k})\ge 2N+2} (-1)^{i(c_k^{2m_k})}
              \dim \overline{C}_{i(c_k^{2m_k})}(E,c_k^{2m_k}; [e])\nn\\
&=& \sum_{k=1}^{q} 2m_k\hat{\chi}(c_k)\nn\\
& & -\sum_{1\le k\le q \atop i(c_k^{2m_k})\ge 2N+2} (-1)^{i(c_k^{2m_k})}
             \dim \overline{C}_{i(c_k^{2m_k})}(E,c_k^{2m_k}; [e]). \lb{pjzbdsbx2}\eea

In order to exactly know whether the iterate $c_k^{2m_k}$ of $c_k$ has contribution to the alternative sum
$\sum_{p=0}^{2N+1}(-1)^p M_p$,  $1\le k\le q$, we set
\bea
N_+^e &=& ^{\#}\{1\le k\le q\ |\ i(c_k^{2m_k})\ge 2N+2,\
i(c_k^{2m_k})-i(c_k)\in 2\Z,\ i(c_k)\in 2\Z\},\lb{jh1}\\
N_+^o &=& ^{\#}\{1\le k\le q\ |\ i(c_k^{2m_k})\ge 2N+2,\
i(c_k^{2m_k})-i(c_k)\in 2\Z,\ i(c_k)\in 2\Z-1\},\lb{jh2}\\
N_-^e &=& ^{\#}\{1\le k\le q\ |\ i(c_k^{2m_k})\le 2N-2,\
i(c_k^{2m_k})-i(c_k)\in 2\Z,\ i(c_k)\in 2\Z\},\lb{jh3}\\
N_-^o &=& ^{\#}\{1\le k\le q\ |\ i(c_k^{2m_k})\le 2N-2,\
i(c_k^{2m_k})-i(c_k)\in 2\Z,\ i(c_k)\in 2\Z-1\}.\lb{jh4}\eea
Here note that $\bar{B}(\Lambda_e M, \Lambda^0 M;\Q)=\frac{n+1}{2n}$,
and  by Theorem \ref{ECIJT}, we can suppose that $N$ is a
multiple of $n$. Thus by Claim 1, (\ref{pjzbdsbx2}), Proposition 2.1,
the definitions of $N^{e}_+$ and $N_+^o$ and Morse inequality, we have
\bea \frac{N(n+1)}{n}+N_+^o-N_+^e
&=& \sum_{k=1}^q 2m_k\hat{\chi}(c_k)-(-N_+^o+N_+^e)  \nn\\
&=& \sum_{p=0}^{2N+1}(-1)^p M_p \nn\\
&\leq& \sum_{p=0}^{2N+1}(-1)^p \beta_p=\frac{N(n+1)}{n}-n,\lb{yibanguji1}\eea
where the Betti numbers $\beta_p$'s are given by (\ref{ksfzbds}). Then
\be N_+^e\ge n.\lb{yibanguji2}\ee

Similar to (\ref{zbgj1})-(\ref{zbgj3}), it follows from Theorem \ref{ECIJT} that there exist also
infinitely many $(q+1)$-tuples $(N', m_1', \cdots, m_q')\in\N^{q+1}$ such that for any $1\le k\le q$,
there holds
\bea
i(c_k^{2m^\prime_k-m}) &\leq& 2N^\prime-2,\quad 1\leq m< 2m^\prime_k, \lb{zbgj1'}\\
i(c_k^{2m^\prime_k}) &=& 2N^\prime-C(M_k)+2\Delta^\prime_k,\lb{zbgj2'}\\
i(c_k^{2m^\prime_k+m})&\geq& 2N^\prime+2,\quad \forall m\geq 1.\lb{zbgj3'}\eea
where, furthermore, $\Delta_k$ and $\Delta_k'$ satisfy the following relationship
\bea \Delta_k' + \Delta_k = C(M_k),\qquad\forall\ 1\le k\le q,\lb{deerta}\eea
which exactly follows from (especially the term (c) of) Theorem 4.2 of \cite{LoZ} (cf. \cite{DLW2} for more details).
Similarly, we define
\bea
N_+^{'e} &=& ^{\#}\{1\le k\le q\ |\ i(c_k^{2m^\prime_k})\ge 2N'+2,\
i(c_k^{2m^\prime_k})-i(c_k)\in 2\Z,\ i(c_k)\in 2\Z\},\lb{jh1'}\\
N_+^{'o} &=& ^{\#}\{1\le k\le q\ |\ i(c_k^{2m^\prime_k})\ge 2N'+2,\
i(c_k^{2m^\prime_k})-i(c_k)\in 2\Z,\ i(c_k)\in 2\Z-1\},\lb{jh2'}\\
N_-^{'e} &=& ^{\#}\{1\le k\le q\ |\ i(c_k^{2m^\prime_k})\le 2N'-2,\
i(c_k^{2m^\prime_k})-i(c_k)\in 2\Z,\ i(c_k)\in 2\Z\},\lb{jh3'}\\
N_-^{'o} &=& ^{\#}\{1\le k\le q\ |\ i(c_k^{2m^\prime_k})\le 2N'-2,\
i(c_k^{2m^\prime_k})-i(c_k)\in 2\Z,\ i(c_k)\in 2\Z-1\}.\lb{jh4'}\eea
So by (\ref{zbgj2'}) and (\ref{deerta}) it yields
\be i(c_k^{2m^\prime_k})= 2N'-C(M_k)+2(C(M_k)-\Delta_k)=2N'+C(M_k)-2\Delta_k. \lb{deertadr}\ee
So by definitions (\ref{jh1})-(\ref{jh4}) and (\ref{jh1'})-(\ref{jh4'}) we have
\be N_{\pm}^e=N_{\mp}^{'e},\qquad N_{\pm}^o=N_{\mp}^{'o}. \lb{dcxz}\ee
Thus, through carrying out the similar arguments to (\ref{yibanguji1})-(\ref{yibanguji2}), by Claim 1,
Proposition 2.1, the definitions of
$N^{'e}_+$ and $N^{'o}_+$  and Morse inequality, we have
\bea \frac{N'(n+1)}{n}+N_+^{'o}-N_+^{'e}
&=& \sum_{k=1}^q 2m'_k\hat{\chi}(c_k)-(-N_+^{'o}+N_+^{'e})  \nn\\
&=& \sum_{p=0}^{2N'+1}(-1)^p M_p \nn\\
&\leq& \sum_{p=0}^{2N'+1}(-1)^p \beta_p=\frac{N'(n+1)}{n}-n,\lb{yibanguji1'}\eea
which, together with (\ref{dcxz}), implies
\be N_-^{e}=N_+^{'e}\ge n.\lb{yibanguji2'}\ee

So by (\ref{yibanguji2}) and (\ref{yibanguji2'}) it yields
\bea q\ge N_+^{e}+N_-^{e}\ge 2n. \lb{n-1}\eea

In addition, any hyperbolic contractible closed geodesic $c_k$ must have $i(c_k^{2m_k})=2N$ since there holds
$C(M_k)=0$ and $\Delta_k=0$ in the hyperbolic case. However, by (\ref{jh1}) and (\ref{jh3}), there exist at least
$(N_+^{e}+N_-^{e})$ contractible closed geodesic with even indices $i(c_k^{2m_k})\geq 2N+2$
or $i(c_k^{2m_k})\leq 2N-2$. So all these
$(N_+^{e}+N_-^{e})$ contractible closed geodesics are non-hyperbolic. Then (\ref{n-1}) shows that there exist
at least $2n$ distinct non-hyperbolic contractible closed geodesics. And (\ref{jh1}), (\ref{jh3}) and (\ref{n-1})
show that all these non-hyperbolic contractible closed geodesics and their iterations have even Morse indices.
This completes the proof of Step 1.

We denote these $2n$ non-hyperbolic contractible closed geodesics by $\{c_k\}_{k=1}^{2n}$.

{\bf Step 2.} {\it There exist at least two distinct contractible closed geodesics different
from those found in Step 1 with even Morse indices.}

\medskip

In fact, for those $2n$ distinct contractible closed geodesics $\{c_k\}_{k=1}^{2n}$ found in Step 1, there holds
$i(c_k^{2m_k})\neq 2N$ by  the definitions of $N_+^e$ and $N_-^e$, which, together with
(\ref{zbgj1}) and (\ref{zbgj3}) yields
\bea i(c_k^m)\neq 2N,\qquad m\ge 1,\quad k=1,\cdots,2n.\lb{bdy2N}\eea
Then by Proposition \ref{Rad1992'} it yields
\bea \sum_{1\le k\le 2n \atop m\ge 1} \dim \overline{C}_{2N}(E, c_k^m; [e])=0.\lb{2Nljwl} \eea
Therefore, noting that $N$ is a multiple of $n$, by (\ref{2Nljwl}), (\ref{ksfzbds})
of Proposition \ref{ksfztd} and Morse inequality, we obtain
\bea \sum_{2n+1\le k\le q}\dim \overline{C}_{2N}(E, c_k^{2m_k}; [e])
&=& \sum_{2n+1\le k\le q,\ m\ge 1}\dim \overline{C}_{2N}(E, c_k^m; [e]) \nn\\
&=& \sum_{1\le k\le q,\ m\ge 1}\dim \overline{C}_{2N}(E, c_k^m; [e])\nn\\
&=& M_{2N}\ge \beta_{2N}=2.\lb{ltbcdx}\eea
where the first equality follows from (\ref{zbgj1}), (\ref{zbgj3}) and Proposition \ref{Rad1992'}.

Now by (\ref{ltbcdx}) and Proposition 2.1, it yields that there exist at least two contractible closed geodesic
$c_{2n+1}$ and $c_{2n+2}$ with $i (c^{2m_k}_k) = 2N$ and
$i (c^{2m_k}_k ) - i (c_k ) \in 2\Z$ for $k = 2n+1$ and $2n+2$. Thus both
$c_{2n+1}$ and $c_{2n+2}$ are different from $\{c_k\}_{k=1}^{2n}$ by (\ref{bdy2N}),
they and all of their iterates have even
Morse indices. This completes the proof of Step 2.

{\bf Acknowledgements.} I would like to thank sincerely the referee for his careful reading of the manuscript, valuable comments on improving the
exposition, and for his deep insight on the main ideas of this paper. And also I would like to sincerely thank Professor Yiming Long,
for introducing me to Hamiltonian dynamics, and for his constant helps and encouragements on my research.

\end{document}